\documentclass[a4paper,12pt]{amsart}

\usepackage{amssymb}
\usepackage{fullpage}
\usepackage[OT4]{fontenc}
\usepackage{texdraw}

\theoremstyle{definition}
\newtheorem{definition}{Definition}
\theoremstyle{plain}
\newtheorem{theorem}[definition]{Theorem}
\newtheorem{conjecture}[definition]{Conjecture}
\newtheorem{corollary}[definition]{Corollary}
\newtheorem{proposition}[definition]{Proposition}
\theoremstyle{remark}
\newtheorem{remark}[definition]{Remark}
\newtheorem{example}[definition]{Example}

\DeclareMathOperator{\edim}{edim}
\DeclareMathOperator{\rank}{rank}
\DeclareMathOperator{\red}{red}
\DeclareMathOperator{\rg}{rg}

\def\adm{\textsc{DimensionBound }}
\def\ass{\longleftarrow}
\def\field{\mathbb{K}}
\def\N{\mathbb{N}}
\def\R{\mathbb{R}}
\def\proj{\mathbb{P}}
\def\Regularity{\textsc{RegularityBound }}
\def\setofdiag{\mathcal{D}}
\def\sspan{\operatorname{span}}
\def\sys{\mathcal{L}}
\def\tto{\longrightarrow}

\newenvironment{algorithm}[1]{\medskip\noindent{\bf Algorithm} #1\\}{\medskip}

\begin{document}

\title{New effective bounds on the dimension\\of a linear system in $\mathbb P^2$}

\author{Marcin Dumnicki, Witold Jarnicki}

\dedicatory{
Institute of Mathematics, Jagiellonian University, \\
Reymonta 4, 30-059 Krak\'ow, Poland \\
}

\thanks{Email address: Witold.Jarnicki@im.uj.edu.pl (corresponding author), Marcin.Dumnicki@im.uj.edu.pl,
phone number: +48-12-6339781, fax number: +48-12-6324372}






\begin{abstract}
The main goal of this paper is to present an algorithm bounding the dimension
of a linear system of curves of given degree (or monomial basis) with
multiple points in general position. As a result we prove
the Hirschowitz--Harbourne Conjecture when the multiplicities of base
points are bounded by $11$.

\medskip

\end{abstract}

\maketitle

\section{Introduction}

Let $\field$ be a field of characteristic zero, $\N=\{0,1,2,\dots\}$, $\N^\ast=\{1,2,3,\dots\}$.
\begin{definition}
Let $D\subset\mathbb N^2$ be finite, let $m_{1},\dots,m_{r} \in \N$, let $p_{1}, \dots, p_{r}
\in \field^{2}$. Define the vector space (over $\field$):
$$
\sys_{D}(m_{1}p_{1},\dots,m_{r}p_{r}) :=
\bigg\{ f = \sum_{\beta \in D} c_{\beta}X^{\beta} \mid c_{\beta} \in \field,
\frac{\partial^{|\alpha|} f}{\partial X^{\alpha}}(p_{j})
= 0, \ |\alpha| < m_{j}, j=1,\dots,r \bigg\}.
$$
Define the \textit{dimension} of the \textit{system of curves
$\sys_{D}(m_{1},\dots,m_{r})$} to be
$$\dim \sys_{D}(m_{1},\dots,m_{r}) := \min_{\{p_{j}\}_{j=1}^{r}, p_{j} \in \field^{2}} \dim_{\field}
\sys_{D}(m_{1}p_{1}, \dots,m_{r}p_{r}) - 1.$$
\end{definition}

\begin{remark}
If points $p_{1},\dots,p_{r}$ are in general position we have
$$\dim \sys_{D}(m_{1},\dots,m_{r}) = \dim_{\field} \sys_{D}(m_{1}p_{1},\dots,m_{r}p_{r}) - 1.$$
\end{remark}

System $\sys_{D}(m_{1},\dots,m_{r})$ can be understood as a vector space
of curves generated by monomials with exponents from $D$ having multiplicities at least
$m_{1},\dots,m_{r}$ in $r$ general points.

Let us assume that $D = \{ \alpha \in \N^{2} \mid |\alpha| \leq d\}$
for some $d \in \N$. Then the space
$\sys_{D}(m_{1},\dots,m_{r})$ can be associated with the linear system over $\proj^{2}$ ($:=\mathbb P^2\mathbb K$)
$$
\sys_{d}(m_{1},\dots,m_{r}) := |dH|(-\sum_{j=1}^{r}m_{j}p_{j}).
$$
This system contains all divisors from the system $|dH|$ ($H$ being a generic
line in $\proj^{2}$) that have multiplicity
at least $m_{j}$ at $p_{j}$, where $p_{1},\dots,p_{r}$ are in general
position. In particular, the space $\sys_{D}(m_{1},\dots,m_{r})$
and the linear system $\sys_{d}(m_{1},\dots,m_{r})$ have equal dimensions.

Intuitively, the dimension of the system $\sys_{D}(m_{1},\dots,m_{r})$
should be equal to the dimension of $V_D:=\sspan \{ X^{\alpha} \mid \alpha \in D\}$ ($=\#D-1$)
minus the number of conditions imposed by the multiplicities $m_{1}, \dots, m_{r}$.
However, the actual dimension may be different.

\begin{definition}
Let $L=\sys_{D}(m_{1},\dots,m_{r})$ be a system of curves. Define
the \textit{expected dimension of $L$}
$$
\edim L := \max \bigg\{ \# D - 1 - \sum_{j=1}^{r} \binom{m_{j}+1}{2}, -1\bigg\}.
$$
\end{definition}

\begin{proposition}\label{linmapping}
For any system of curves $L$ we have
$$
\dim L \geq \edim L.
$$
\end{proposition}

\begin{proof}
Fix $p_1,\dots,p_r\in\field^2$. Put
$$
\varphi_{j,\alpha}^{p_1,\dots,p_r}:V_D\ni f\mapsto \frac{\partial^{|\alpha|} f}{\partial X^{\alpha}}(p_{j})\in\field,\quad j=1,\dots,r,\;\alpha\in\N^2.
$$
Let $\mathcal A=\{(j,\alpha)\mid|\alpha|<m_j,\;j=1,\dots,r\}$.
Consider the linear mapping
$$
\varphi^{p_1,\dots,p_r}:V_D\ni f\mapsto(\varphi_a^{p_1,\dots,p_r}(f))_{a\in\mathcal A}\in\field^{\#\mathcal A}.
$$
By linear algebra we have $\dim\ker\varphi^{p_1,\dots,p_r}\geq\#D-\#\mathcal A$, which implies
$$
\dim L=\min_{\{p_{j}\}_{j=1}^{r}, p_{j} \in \field^{2}}\dim\ker\varphi^{p_1,\dots,p_r}-1\geq\edim L.
$$
 \end{proof}

\begin{definition}
We say that \textit{system of curves $L$ is special}
if
$$\dim L > \edim L.$$
Otherwise we say that \textit{$L$ is non-special}.
\end{definition}

\section{The Hirschowitz--Harbourne Conjecture}

For systems of the form $\sys_{d}(m_{1},\dots,m_{r})$ the well-known
Hirschowitz--Harbourne Conjecture giving geometrical
description to the speciality of a system  was formulated in \cite{HCON}.
To formulate this Conjecture consider the blowing-up
$\pi : \widetilde\proj^{2} \to \proj$ in $r$ general points with
exceptional divisors $E_{1},\dots,E_{r}$.

\begin{definition}
A curve $C \subset \proj^{2}$ is said to be \textit{$-1$-curve}
if it is irreducible, and the self-intersection of its proper
transform $\widetilde{C} \subset \widetilde\proj^{2}$ is equal to $-1$.
\end{definition}

\begin{conjecture}[Hirschowitz--Harbourne]\label{hirsch}
A system $L = \sys_{d}(m_{1},\dots,m_{r})$ is special if and only if
there exists a $-1$-curve $C \subset \proj^{2}$ such that
$$\widetilde{L}.\widetilde{C} \leq -2,$$
where $\widetilde{L} = |d\pi^{*}(\mathcal O_{\mathbb P^2}(1)) - \sum_{j=1}^{r}m_{j}E_{j}|$.
\end{conjecture}

This Conjecture was studied by many authors, we refer only to the recent
results. For homogenous systems ($m_1=m_{2}=m_{3}=\dots=m_{r}=:m$), the above conjecture holds for $m\leq20$ (see \cite{CCMO,CMir}). In the general case multiplicities bounded by $7$ have been dealt with (see \cite{Yang}). Both these results were obtained with the help of computers.

For further information about above conjecture see for example
\cite{CIL}.

The main result of this paper is Corollary \ref{hirsch11} stating that Conjecture \ref{hirsch} holds for multiplicities bounded by $11$.

For a given system one can verify whether the above conjecture holds. If a multiple $-1$-curve in the base locus cannot be easily found computationally, we check whether the system is non-special by using our reduction method and/or calculating the determinant of $\varphi$ (see Proposition \ref{linmapping}).

Hence, in order to prove Conjecture \ref{hirsch} for linear systems up to multiplicity $11$, it is enough to limit the number of cases to be checked. Therefore we will not focus on the conjecture itself, but on the possibility of
finding the list of all special systems (in particular we use the reduction method to obtain effective bounds on the number of cases).

\section{Reduction method}

\begin{definition}
Let $a_{1}, \dots, a_{k} \in \N$, $a_{j} \leq j$, $j=1,\dots,k$.
Define the diagram $(a_{1},\dots,a_{k})$
$$(a_{1},\dots,a_{k}) = \bigcup_{j=1}^{k} \{ (\alpha_{1},\alpha_{2}) \in \N^{2} \ | \ \alpha_{1} + \alpha_{2} = j - 1, \
\alpha_{2} < a_{j} \}.$$
For $a, a_{1}, \dots, a_{k} \in \N$, $a_{j} \leq a+j$ define
$$(\overline{a},a_{1},\dots,a_{k}) := (1,2,\dots,a-1,a,a_{1},\dots,a_{k}).$$
\end{definition}

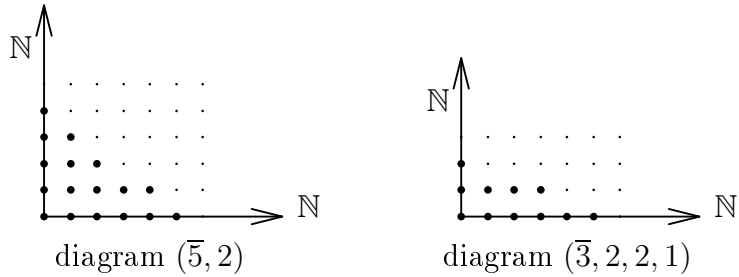
\begin{example}
$$
\begin{array}{ccc}


\begin{texdraw}
\drawdim pt
\move(0 0)
\fcir f:0 r:0.5
\move(0 10)
\fcir f:0 r:0.5
\move(0 20)
\fcir f:0 r:0.5
\move(0 30)
\fcir f:0 r:0.5
\move(0 40)
\fcir f:0 r:0.5
\move(0 50)
\fcir f:0 r:0.5
\move(10 0)
\fcir f:0 r:0.5
\move(10 10)
\fcir f:0 r:0.5
\move(10 20)
\fcir f:0 r:0.5
\move(10 30)
\fcir f:0 r:0.5
\move(10 40)
\fcir f:0 r:0.5
\move(10 50)
\fcir f:0 r:0.5
\move(20 0)
\fcir f:0 r:0.5
\move(20 10)
\fcir f:0 r:0.5
\move(20 20)
\fcir f:0 r:0.5
\move(20 30)
\fcir f:0 r:0.5
\move(20 40)
\fcir f:0 r:0.5
\move(20 50)
\fcir f:0 r:0.5
\move(30 0)
\fcir f:0 r:0.5
\move(30 10)
\fcir f:0 r:0.5
\move(30 20)
\fcir f:0 r:0.5
\move(30 30)
\fcir f:0 r:0.5
\move(30 40)
\fcir f:0 r:0.5
\move(30 50)
\fcir f:0 r:0.5
\move(40 0)
\fcir f:0 r:0.5
\move(40 10)
\fcir f:0 r:0.5
\move(40 20)
\fcir f:0 r:0.5
\move(40 30)
\fcir f:0 r:0.5
\move(40 40)
\fcir f:0 r:0.5
\move(40 50)
\fcir f:0 r:0.5
\move(50 0)
\fcir f:0 r:0.5
\move(50 10)
\fcir f:0 r:0.5
\move(50 20)
\fcir f:0 r:0.5
\move(50 30)
\fcir f:0 r:0.5
\move(50 40)
\fcir f:0 r:0.5
\move(50 50)
\fcir f:0 r:0.5
\move(60 0)
\fcir f:0 r:0.5
\move(60 10)
\fcir f:0 r:0.5
\move(60 20)
\fcir f:0 r:0.5
\move(60 30)
\fcir f:0 r:0.5
\move(60 40)
\fcir f:0 r:0.5
\move(60 50)
\fcir f:0 r:0.5
\arrowheadtype t:V
\move(0 0)
\avec(90 0)
\move(0 0)
\avec(0 80)
\htext(96 0){$\mathbb{N}$}
\htext(-13 60){$\mathbb{N}$}
\move(0 0)
\fcir f:0 r:1.5
\move(10 0)
\fcir f:0 r:1.5
\move(0 10)
\fcir f:0 r:1.5
\move(20 0)
\fcir f:0 r:1.5
\move(10 10)
\fcir f:0 r:1.5
\move(0 20)
\fcir f:0 r:1.5
\move(30 0)
\fcir f:0 r:1.5
\move(20 10)
\fcir f:0 r:1.5
\move(10 20)
\fcir f:0 r:1.5
\move(0 30)
\fcir f:0 r:1.5
\move(40 0)
\fcir f:0 r:1.5
\move(30 10)
\fcir f:0 r:1.5
\move(20 20)
\fcir f:0 r:1.5
\move(10 30)
\fcir f:0 r:1.5
\move(0 40)
\fcir f:0 r:1.5
\move(50 0)
\fcir f:0 r:1.5
\move(40 10)
\fcir f:0 r:1.5
\textref h:C v:C
\htext(40 -15){diagram $(\overline{5},2)$}
\end{texdraw}


&
\hspace{1cm}
&


\begin{texdraw}
\drawdim pt
\move(0 0)
\fcir f:0 r:0.5
\move(0 10)
\fcir f:0 r:0.5
\move(0 20)
\fcir f:0 r:0.5
\move(0 30)
\fcir f:0 r:0.5
\move(10 0)
\fcir f:0 r:0.5
\move(10 10)
\fcir f:0 r:0.5
\move(10 20)
\fcir f:0 r:0.5
\move(10 30)
\fcir f:0 r:0.5
\move(20 0)
\fcir f:0 r:0.5
\move(20 10)
\fcir f:0 r:0.5
\move(20 20)
\fcir f:0 r:0.5
\move(20 30)
\fcir f:0 r:0.5
\move(30 0)
\fcir f:0 r:0.5
\move(30 10)
\fcir f:0 r:0.5
\move(30 20)
\fcir f:0 r:0.5
\move(30 30)
\fcir f:0 r:0.5
\move(40 0)
\fcir f:0 r:0.5
\move(40 10)
\fcir f:0 r:0.5
\move(40 20)
\fcir f:0 r:0.5
\move(40 30)
\fcir f:0 r:0.5
\move(50 0)
\fcir f:0 r:0.5
\move(50 10)
\fcir f:0 r:0.5
\move(50 20)
\fcir f:0 r:0.5
\move(50 30)
\fcir f:0 r:0.5
\move(60 0)
\fcir f:0 r:0.5
\move(60 10)
\fcir f:0 r:0.5
\move(60 20)
\fcir f:0 r:0.5
\move(60 30)
\fcir f:0 r:0.5
\arrowheadtype t:V
\move(0 0)
\avec(90 0)
\move(0 0)
\avec(0 60)
\htext(96 0){$\mathbb{N}$}
\htext(-13 40){$\mathbb{N}$}
\move(0 0)
\fcir f:0 r:1.5
\move(10 0)
\fcir f:0 r:1.5
\move(0 10)
\fcir f:0 r:1.5
\move(20 0)
\fcir f:0 r:1.5
\move(10 10)
\fcir f:0 r:1.5
\move(0 20)
\fcir f:0 r:1.5
\move(30 0)
\fcir f:0 r:1.5
\move(20 10)
\fcir f:0 r:1.5
\move(40 0)
\fcir f:0 r:1.5
\move(30 10)
\fcir f:0 r:1.5
\move(50 0)
\fcir f:0 r:1.5
\textref h:C v:C
\htext(40 -15){diagram $(\overline{3},2,2,1)$}
\end{texdraw}


\end{array}
$$
\end{example}

\begin{definition}
Let $k \in \N^{*}$, let $D = (b_{1},\dots,b_{\ell},a_{1},\dots,a_{k})$ be a diagram, $a_{k} > 0$.
Define the numbers $v_{j}$, $j=k,\dots,1$ inductively to be
$$v_{j} := \left\{
\begin{array}{ll}
a_{j}, & a_{j} < k, \\
\max \{ 1,\dots,k \} \setminus \{v_{j+1},\dots,v_{k}\},& a_{j} \geq k.
\end{array}
\right.$$
If we have
$$v_{i} \neq v_{j} \textrm{ for }i \neq j, \qquad v_{j} \leq a_{j}
\textrm{ for }j=1,\dots,k$$
then we say that
\textit{$D$ is $k$-reducible}. The diagram
$$\red_{k}(D) := (b_{1},\dots,b_{\ell},
a_{1}-v_{1}, \dots, a_{k}-v_{k})$$
will be called the \textit{$k$-reduction of $D$}. We use the following notation
$$
(b_1,\dots,b_\ell,a_1,\dots,a_k)\stackrel k\tto(b_1,\dots,b_\ell,a_1-v_1,\dots,a_k-v_k).
$$
\end{definition}

\begin{definition}\label{weakred}
Let $k \in \N^{*}$.
We say that a diagram $G=(b_{1},\dots,b_{\ell},a_{1},\dots,a_{k})$ is
\textit{a weak $k$-reduction of a diagram
$D=(b_{1},\dots,b_{\ell},a_{1}',\dots,a_{k}')$}
if there exists a sequence of numbers $c_{1}, \dots,c_{k} \in \N$ (possibly all
of them equal to $0$) such that
$$
D':=(b_{1},\dots,b_{\ell},a_{1}'+c_{1},\dots,
a_{k}'+c_{k})
$$
is a diagram (in particular $a_j'+c_j\leq\ell+j$, $j=1,\dots, k$) and
$\red_{k}(D')=G$. We use the following notation
$$
(b_1,\dots,b_\ell,a_1',\dots,a_k')\stackrel{k\textrm{ w}}\tto(b_1,\dots,b_\ell,a_1,\dots,a_k).
$$
\end{definition}

\begin{remark}Observe that the $k$-reduction of $D$ is a weak $k$-reduction.
\end{remark}

\begin{example}
$$
\begin{array}{ccc}


\begin{texdraw}
\drawdim pt
\move(0 0)
\fcir f:0 r:0.5
\move(0 10)
\fcir f:0 r:0.5
\move(0 20)
\fcir f:0 r:0.5
\move(0 30)
\fcir f:0 r:0.5
\move(0 40)
\fcir f:0 r:0.5
\move(0 50)
\fcir f:0 r:0.5
\move(10 0)
\fcir f:0 r:0.5
\move(10 10)
\fcir f:0 r:0.5
\move(10 20)
\fcir f:0 r:0.5
\move(10 30)
\fcir f:0 r:0.5
\move(10 40)
\fcir f:0 r:0.5
\move(10 50)
\fcir f:0 r:0.5
\move(20 0)
\fcir f:0 r:0.5
\move(20 10)
\fcir f:0 r:0.5
\move(20 20)
\fcir f:0 r:0.5
\move(20 30)
\fcir f:0 r:0.5
\move(20 40)
\fcir f:0 r:0.5
\move(20 50)
\fcir f:0 r:0.5
\move(30 0)
\fcir f:0 r:0.5
\move(30 10)
\fcir f:0 r:0.5
\move(30 20)
\fcir f:0 r:0.5
\move(30 30)
\fcir f:0 r:0.5
\move(30 40)
\fcir f:0 r:0.5
\move(30 50)
\fcir f:0 r:0.5
\move(40 0)
\fcir f:0 r:0.5
\move(40 10)
\fcir f:0 r:0.5
\move(40 20)
\fcir f:0 r:0.5
\move(40 30)
\fcir f:0 r:0.5
\move(40 40)
\fcir f:0 r:0.5
\move(40 50)
\fcir f:0 r:0.5
\move(50 0)
\fcir f:0 r:0.5
\move(50 10)
\fcir f:0 r:0.5
\move(50 20)
\fcir f:0 r:0.5
\move(50 30)
\fcir f:0 r:0.5
\move(50 40)
\fcir f:0 r:0.5
\move(50 50)
\fcir f:0 r:0.5
\move(60 0)
\fcir f:0 r:0.5
\move(60 10)
\fcir f:0 r:0.5
\move(60 20)
\fcir f:0 r:0.5
\move(60 30)
\fcir f:0 r:0.5
\move(60 40)
\fcir f:0 r:0.5
\move(60 50)
\fcir f:0 r:0.5
\move(70 0)
\fcir f:0 r:0.5
\move(70 10)
\fcir f:0 r:0.5
\move(70 20)
\fcir f:0 r:0.5
\move(70 30)
\fcir f:0 r:0.5
\move(70 40)
\fcir f:0 r:0.5
\move(70 50)
\fcir f:0 r:0.5
\arrowheadtype t:V
\move(0 0)
\avec(100 0)
\move(0 0)
\avec(0 80)
\htext(106 0){$\mathbb{N}$}
\htext(-13 60){$\mathbb{N}$}
\move(0 0)
\fcir f:0 r:1.5
\move(10 0)
\fcir f:0 r:1.5
\move(0 10)
\fcir f:0 r:1.5
\move(20 0)
\fcir f:0 r:1.5
\move(10 10)
\fcir f:0 r:1.5
\move(0 20)
\fcir f:0 r:1.5
\move(30 0)
\fcir f:0 r:1.5
\move(20 10)
\fcir f:0 r:1.5
\move(10 20)
\fcir f:0 r:1.5
\move(0 30)
\fcir f:0 r:1.5
\move(40 0)
\fcir f:0 r:1.5
\move(30 10)
\fcir f:0 r:1.5
\move(20 20)
\fcir f:0 r:1.5
\move(10 30)
\fcir f:0 r:1.5
\move(0 40)
\fcir f:0 r:1.5
\move(50 0)
\fcir f:0 r:1.5
\move(40 10)
\fcir f:0 r:1.5
\move(30 20)
\fcir f:0 r:1.5
\move(60 0)
\fcir f:0 r:1.5
\move(50 10)
\fcir f:0 r:1.5
\move(0 30)
\move(-3 27)
\lvec(3 33)
\move(3 27)
\lvec(-3 33)
\move(0 40)
\move(-3 37)
\lvec(3 43)
\move(3 37)
\lvec(-3 43)
\move(10 30)
\move(7 27)
\lvec(13 33)
\move(13 27)
\lvec(7 33)
\move(20 20)
\move(17 17)
\lvec(23 23)
\move(23 17)
\lvec(17 23)
\move(30 10)
\move(27 7)
\lvec(33 13)
\move(33 7)
\lvec(27 13)
\move(30 20)
\move(27 17)
\lvec(33 23)
\move(33 17)
\lvec(27 23)
\move(40 10)
\move(37 7)
\lvec(43 13)
\move(43 7)
\lvec(37 13)
\move(50 0)
\move(47 -3)
\lvec(53 3)
\move(53 -3)
\lvec(47 3)
\move(50 10)
\move(47 7)
\lvec(53 13)
\move(53 7)
\lvec(47 13)
\move(60 0)
\move(57 -3)
\lvec(63 3)
\move(63 -3)
\lvec(57 3)
\textref h:C v:C
\htext(45 -15){the $4$-reduction of diagram $(\overline{5},3,2)$}
\end{texdraw}


&
\hspace{2cm}
&


\begin{texdraw}
\drawdim pt
\move(0 0)
\fcir f:0 r:0.5
\move(0 10)
\fcir f:0 r:0.5
\move(0 20)
\fcir f:0 r:0.5
\move(0 30)
\fcir f:0 r:0.5
\move(0 40)
\fcir f:0 r:0.5
\move(10 0)
\fcir f:0 r:0.5
\move(10 10)
\fcir f:0 r:0.5
\move(10 20)
\fcir f:0 r:0.5
\move(10 30)
\fcir f:0 r:0.5
\move(10 40)
\fcir f:0 r:0.5
\move(20 0)
\fcir f:0 r:0.5
\move(20 10)
\fcir f:0 r:0.5
\move(20 20)
\fcir f:0 r:0.5
\move(20 30)
\fcir f:0 r:0.5
\move(20 40)
\fcir f:0 r:0.5
\move(30 0)
\fcir f:0 r:0.5
\move(30 10)
\fcir f:0 r:0.5
\move(30 20)
\fcir f:0 r:0.5
\move(30 30)
\fcir f:0 r:0.5
\move(30 40)
\fcir f:0 r:0.5
\move(40 0)
\fcir f:0 r:0.5
\move(40 10)
\fcir f:0 r:0.5
\move(40 20)
\fcir f:0 r:0.5
\move(40 30)
\fcir f:0 r:0.5
\move(40 40)
\fcir f:0 r:0.5
\move(50 0)
\fcir f:0 r:0.5
\move(50 10)
\fcir f:0 r:0.5
\move(50 20)
\fcir f:0 r:0.5
\move(50 30)
\fcir f:0 r:0.5
\move(50 40)
\fcir f:0 r:0.5
\move(60 0)
\fcir f:0 r:0.5
\move(60 10)
\fcir f:0 r:0.5
\move(60 20)
\fcir f:0 r:0.5
\move(60 30)
\fcir f:0 r:0.5
\move(60 40)
\fcir f:0 r:0.5
\arrowheadtype t:V
\move(0 0)
\avec(90 0)
\move(0 0)
\avec(0 70)
\htext(96 0){$\mathbb{N}$}
\htext(-13 50){$\mathbb{N}$}
\move(0 0)
\fcir f:0 r:1.5
\move(10 0)
\fcir f:0 r:1.5
\move(0 10)
\fcir f:0 r:1.5
\move(20 0)
\fcir f:0 r:1.5
\move(10 10)
\fcir f:0 r:1.5
\move(0 20)
\fcir f:0 r:1.5
\move(30 0)
\fcir f:0 r:1.5
\move(20 10)
\fcir f:0 r:1.5
\move(10 20)
\fcir f:0 r:1.5
\move(0 30)
\fcir f:0 r:1.5
\move(40 0)
\fcir f:0 r:1.5
\move(30 10)
\fcir f:0 r:1.5
\move(20 20)
\fcir f:0 r:1.5
\move(10 30)
\fcir f:0 r:1.5
\move(50 0)
\fcir f:0 r:1.5
\move(40 10)
\fcir f:0 r:1.5
\move(30 20)
\fcir f:0 r:1.5
\move(20 30)
\fcir f:0 r:1.5
\move(0 10)
\move(-3 7)
\lvec(3 13)
\move(3 7)
\lvec(-3 13)
\move(0 20)
\move(-3 17)
\lvec(3 23)
\move(3 17)
\lvec(-3 23)
\move(10 10)
\move(7 7)
\lvec(13 13)
\move(13 7)
\lvec(7 13)
\move(0 30)
\move(-3 27)
\lvec(3 33)
\move(3 27)
\lvec(-3 33)
\move(10 20)
\move(7 17)
\lvec(13 23)
\move(13 17)
\lvec(7 23)
\move(20 10)
\move(17 7)
\lvec(23 13)
\move(23 7)
\lvec(17 13)
\move(10 30)
\move(7 27)
\lvec(13 33)
\move(13 27)
\lvec(7 33)
\move(20 20)
\move(17 17)
\lvec(23 23)
\move(23 17)
\lvec(17 23)
\move(30 10)
\move(27 7)
\lvec(33 13)
\move(33 7)
\lvec(27 13)
\move(40 0)
\move(37 -3)
\lvec(43 3)
\move(43 -3)
\lvec(37 3)
\move(20 30)
\move(17 27)
\lvec(23 33)
\move(23 27)
\lvec(17 33)
\move(30 20)
\move(27 17)
\lvec(33 23)
\move(33 17)
\lvec(27 23)
\move(40 10)
\move(37 7)
\lvec(43 13)
\move(43 7)
\lvec(37 13)
\move(50 0)
\move(47 -3)
\lvec(53 3)
\move(53 -3)
\lvec(47 3)
\textref h:C v:C
\htext(60 -15){a weak $5$-reduction of diagram $(\overline{4},4,4)$}
\end{texdraw}


\end{array}
$$
\end{example}

Now we can formulate the main theorem.

\begin{theorem}
\label{mainred}
Let $m_{1},\dots,m_{r} \in \N$. If a diagram $D'$ is a
weak $m_{r}$-reduction of a diagram $D$ then
$$\dim \sys_{D}(m_{1},\dots,m_{r}) \leq \dim \sys_{D'}(m_{1},\dots,m_{r-1}).$$
\end{theorem}

\begin{definition}
In the situation above, we say that \textit{$\sys_{D'}(m_{1},\dots,m_{r-1})$ is
a weak reduction of the system $\sys_{D}(m_{1},\dots,m_{r})$}, or simply
\textit{the reduction} if $D'=\red_{m_{r}}(D)$.
\end{definition}

In particular, we have

\begin{corollary}
Let $m_{1},\dots,m_{r} \in \N$. Let $L=\sys_{D}(m_{1},\dots,m_{r})$.
If $D$ is $m_{r}$-reducible and the reduction of $L$ is non-special, then $L$ is non-special.
\end{corollary}

\begin{proof}[{\it Proof of Theorem \ref{mainred}}]
We present here a sketch of the proof (for details see
\cite{Dum}).
If $D'$ is a weak $m_{r}$-reduction of $D$ then there exists
a diagram $G$ such that $D \subset G$, and $D'$ is the $m_{r}$-reduction
of $G$. Of course
$$\dim \sys_{D}(m_{1},\dots,m_{r}) \leq \dim \sys_{G}(m_{1},\dots,m_{r}).$$
It remains to prove the inequality
$$\dim \sys_{G}(m_{1},\dots,m_{r}) \leq \dim \sys_{D'}(m_{1},\dots,m_{r-1}).$$
We proceed in two steps.\\
\textit{Step 1.}
We will show the following:
\begin{enumerate}
\item
\label{p1}
the system $\sys_{G \setminus D'}(m_{r})$ is non-special,
\item
\label{p2}
if a set $P \subset G$ satisfies $\# P = \# (G\setminus D')$ and
$$\sum_{\alpha \in P} \alpha = \sum_{\alpha \in
G \setminus D'} \alpha$$
then the system $\sys_{P}(m_{r})$ is special.
\end{enumerate}
To prove (1), observe that the system $\sys_{G \setminus D'}(m_{r})$ is non-special
if and only if the elements of $G \setminus D'$, considered as points in $\N^{2} \subset
\R^{2}$, do not lie on a curve of degree $m_{r} - 1$. The last condition
holds in view of Bezout's Theorem.
To prove (2), consider the set
$$
\mathcal F= \bigg\{ F \subset G\mid \#F = \binom{m_{r}+1}{2},
F \textrm{ do not lie on a curve of degree } m_{r}-1\bigg\}.
$$
For $F\in\mathcal F$ put $|F|:=\sum_{\alpha \in F}|\alpha|$.
It can be shown by induction on $m_{r}$ that $K:=G \setminus D'$ is the only element in $\mathcal F$ such that $|K|\geq|F|$ for all $F\in\mathcal F$ (this is the most technical part of the proof, for details see \cite{Dum}).

\textit{Step 2.}
Now we can compute the dimension of a system $L:=\sys_{G}(m_{1},\dots,m_{r})$ as
$\dim L = \# G - \rank M - 1$, where $M$ is the matrix of the mapping $\varphi^{p_1,\dots,p_r}$ associated to $L$ (see the proof of Proposition \ref{linmapping}). In an appropriate basis this matrix is of the following form
$$M = \left[ \begin{array}{c|c}
\textrm{ matrix } & \\
\textrm{ associated to } & K \\
\sys_{D'}(m_{1}p_{1},\dots,m_{r-1}p_{r-1}) & \\ \hline
& \textrm{ square matrix } \\
K' & \textrm{ associated to } \\
& \sys_{G \setminus D'}(m_{r}p_{r}) \\ \end{array} \right].$$
The determinant of the lower right submatrix is nonzero (because of (\ref{p1})).
Take a maximal nonzero minor (of size $k$) in the upper left submatrix.
Then, by generalized Laplace rule and property (\ref{p2}),
$$
\rank M = k + \# G - \# D',
$$
which finishes the proof.
\end{proof}

In the proof we have used two properties of the set of reduced monomials.
In fact these two properties may hold for many other sets, which
allows us to find and use other ``reductions''.

\begin{definition}
Let $D\subset\N^{2}$ be finite (not necessarily a diagram).
We say that a system $\sys_{D}(m_{1},\dots,m_{r})$
\textit{admits a reduction algorithm} if there exists
a sequence of sets $D = D_{r} \supset \dots \supset D_{0}$ such that
\begin{enumerate}
\item $\# D_{j} - \# D_{j-1} \leq \binom{m_{j}+1}{2}$, $j=1,\dots,r$,
\item the system $\sys_{D_{j} \setminus D_{j-1}}(m_{j})$ is non-special, $j=1,\dots,r$,
\item
if $P \subset D_{j}$ satisfies $\# P = \# D_{j} - \# D_{j-1}$ and
$$\sum_{\alpha \in P} \alpha
= \sum_{\alpha \in D_{j} \setminus D_{j-1}} \alpha$$
then $\sys_{P}(m_{j})$ is special, $j=1,\dots,r$,
\item $\# D_{0} - 1 \leq \edim \sys_{D}(m_{1},\dots,m_{r})$.
\end{enumerate}
\end{definition}

The following theorem has been proven in
\cite{Dum} (also for dimension greater than two).

\begin{theorem}
\label{algnos}
If a system admits a reduction algorithm, then it is non-special.
\end{theorem}

Moreover, investigating non-speciality of many systems (in dimension two, as well as for higher dimensions) has lead to the following

\begin{conjecture}
Every non-special system admits a reduction algorithm.
\end{conjecture}

\begin{remark}
In
\cite{Yang}
S. Yang presented the ``box diagram algorithm'', which can
be treated as another possibility of reducing diagrams. However, not
every non-special system can be reduced to a trivial one using
the reduction of this type. Moreover, Theorem \ref{algnos} works
also in a higher dimension.
\end{remark}

Now we can present the main algorithm to bound the dimension of a system.

\begin{algorithm}{\adm}

\noindent {\bf Input:} a diagram $D = (\overline{a},a_{1},\dots,a_{k})$
and sequence of multiplicities $(m_{1},\dots,m_{r})$. \\
{\bf Output:} a number $d \in \N$ such that
$\dim \sys_{D}(m_{1},\dots,m_{r}) \leq d$. \\

\medskip

\noindent {\bf for} $i \ass 1$ {\bf to} $r$ {\bf do $\{$}

\nopagebreak
$D \ass$ a weak $m_{i}$-reduction of $D$ \\
$\}$ \\
$d \ass \# D - 1$ \\
{\bf return} $d$

\end{algorithm}

\begin{remark}Observe that the quality of the bound obtained depends on the weak reductions chosen. In particular taking $D$ as a weak reduction of $D$ results in a trivial bound on the dimension. Our implementation of this algorithm chooses $D'\supset D$ (see Definition \ref{weakred}) of minimum possible cardinality, which results in much better bounds.
\end{remark}

\begin{example}
\label{ex27}
Consider the system $L = \sys_{27}(24,6,6,6,6,6,6)$. We can see
that $\dim L \geq 0$ (picture). The following sequence of
weak reductions
\begin{align*}
(\overline{28}) &\stackrel{24}{\tto}
(\overline{4},\underbrace{4,\dots,4}_{24}) \stackrel{6\textrm{ w}}{\tto}
(\overline{4},\underbrace{4,\dots,4}_{18},3,2,1) \stackrel{6\textrm{ w}}{\tto}
(\overline{4},\underbrace{4,\dots,4}_{15})
\\
&\stackrel{6\textrm{ w}}{\tto}
(\overline{4},\underbrace{4,\dots,4}_{9},3,2,1) \stackrel{6\textrm{ w}}{\tto} (\overline{4},4,4,4,4,4,4) \stackrel{6\textrm{ w}}{\tto}
(\overline{4},3,2,1) \stackrel{6\textrm{ w}}{\tto}
(1)
\end{align*}
gives $\dim L \leq 0$.

\medskip
\centertexdraw{
\drawdim pt
\linewd 0.8
\move(35.0 0.0)
\clvec(35.66 1.66)(36.33 3.0)(37.33 4.33)
\move(37.33 4.33)
\clvec(37.66 6.0)(38.66 7.33)(39.33 9.0)
\move(39.33 9.0)
\clvec(40.33 10.33)(40.66 12.0)(42.0 13.33)
\move(42.0 13.33)
\clvec(42.33 15.0)(43.33 16.33)(44.0 18.0)
\move(44.0 18.0)
\clvec(44.66 19.33)(45.33 20.66)(46.0 22.33)
\move(46.0 22.33)
\clvec(47.0 23.66)(47.33 25.33)(48.33 26.66)
\move(48.33 26.66)
\clvec(48.66 28.33)(49.66 29.66)(50.0 31.33)
\move(50.0 31.33)
\clvec(50.66 32.66)(51.0 34.33)(52.0 35.66)
\move(52.0 35.66)
\clvec(52.33 37.33)(53.0 39.0)(53.33 40.66)
\move(53.33 40.66)
\clvec(53.66 42.33)(54.0 44.0)(54.33 45.66)
\move(54.33 45.66)
\clvec(54.33 47.33)(55.0 49.0)(54.66 50.66)
\move(54.66 50.66)
\clvec(54.33 52.33)(54.0 54.0)(53.66 55.66)
\move(53.66 55.66)
\clvec(53.0 57.33)(52.33 59.0)(51.66 60.66)
\move(51.66 60.66)
\clvec(51.0 62.0)(50.33 63.33)(49.33 64.66)
\move(49.33 64.66)
\clvec(48.33 66.0)(47.33 67.33)(46.33 68.66)
\move(46.33 68.66)
\clvec(45.66 70.33)(44.66 71.66)(43.66 73.0)
\move(43.66 73.0)
\clvec(42.66 74.33)(41.66 75.66)(40.66 77.0)
\move(40.66 77.0)
\clvec(40.0 78.33)(39.0 79.66)(38.0 81.0)
\move(38.0 81.0)
\clvec(37.33 82.66)(36.33 84.0)(35.66 85.66)
\move(35.66 85.66)
\clvec(34.66 87.0)(34.33 88.66)(33.33 90.0)
\move(33.33 90.0)
\clvec(33.0 91.66)(32.33 93.0)(32.0 94.66)
\move(32.0 94.66)
\clvec(31.66 96.33)(31.33 98.0)(31.33 99.66)
\move(31.33 99.66)
\clvec(31.33 101.33)(31.66 103.0)(32.0 104.66)
\move(32.0 104.66)
\clvec(32.33 106.33)(33.0 107.66)(34.0 109.0)
\move(34.0 109.0)
\clvec(35.0 110.33)(36.33 111.66)(37.66 112.66)
\move(37.66 112.66)
\clvec(39.0 113.33)(40.33 114.66)(42.0 115.0)
\move(42.0 115.0)
\clvec(43.33 115.66)(45.0 116.0)(46.66 116.33)
\move(46.66 116.33)
\clvec(48.33 116.66)(50.0 117.0)(51.66 117.33)
\move(51.66 117.33)
\clvec(53.33 117.33)(55.0 118.0)(56.66 117.66)
\move(56.66 117.66)
\clvec(58.33 118.0)(60.0 117.66)(61.66 118.0)
\move(61.66 118.0)
\clvec(63.33 117.66)(65.0 117.66)(66.66 117.33)
\move(66.66 117.33)
\clvec(68.33 117.33)(70.0 117.0)(71.66 116.66)
\move(71.66 116.66)
\clvec(73.33 116.33)(75.0 116.0)(76.66 115.66)
\move(76.66 115.66)
\clvec(78.33 115.33)(80.0 115.0)(81.66 114.66)
\move(81.66 114.66)
\clvec(83.33 114.33)(85.0 114.0)(86.66 113.33)
\move(86.66 113.33)
\clvec(88.33 113.33)(90.0 112.66)(91.66 112.0)
\move(91.66 112.0)
\clvec(93.33 111.66)(95.0 111.0)(96.66 110.33)
\move(96.66 110.33)
\clvec(98.33 110.0)(100.0 109.33)(101.66 109.0)
\move(101.66 109.0)
\clvec(103.0 108.0)(104.66 107.66)(106.0 107.0)
\move(106.0 107.0)
\clvec(107.66 106.33)(109.33 105.66)(110.66 105.0)
\move(110.66 105.0)
\clvec(112.0 104.33)(113.66 103.66)(115.0 103.0)
\move(115.0 103.0)
\clvec(116.33 102.33)(118.0 101.66)(119.33 100.66)
\move(119.33 100.66)
\clvec(121.0 100.33)(122.33 99.33)(123.66 98.66)
\move(123.66 98.66)
\clvec(125.0 97.66)(126.66 97.0)(128.0 96.0)
\move(128.0 96.0)
\clvec(129.66 95.33)(131.0 94.33)(132.33 93.33)
\move(132.33 93.33)
\clvec(134.0 92.66)(135.33 91.66)(136.66 90.66)
\move(136.66 90.66)
\clvec(138.0 89.66)(139.33 89.0)(140.66 88.0)
\move(140.66 88.0)
\clvec(142.0 87.0)(143.33 86.0)(144.66 85.0)
\move(144.66 85.0)
\clvec(146.0 84.0)(147.33 83.0)(148.66 81.66)
\move(148.66 81.66)
\clvec(150.0 81.0)(151.33 79.66)(152.66 78.66)
\move(152.66 78.66)
\clvec(153.66 77.33)(155.0 76.33)(156.33 75.0)
\move(156.33 75.0)
\clvec(157.33 73.66)(158.66 72.66)(159.66 71.33)
\move(159.66 71.33)
\clvec(160.66 70.0)(161.66 68.66)(163.0 67.33)
\move(163.0 67.33)
\clvec(164.0 66.0)(164.66 64.66)(166.0 63.33)
\move(166.0 63.33)
\clvec(166.66 61.66)(167.66 60.33)(168.0 58.66)
\move(168.0 58.66)
\clvec(169.0 57.33)(169.33 55.66)(169.66 54.0)
\move(169.66 54.0)
\clvec(170.0 52.33)(170.33 50.66)(170.33 49.0)
\move(170.33 49.0)
\clvec(170.0 47.33)(169.66 45.66)(169.0 44.0)
\move(169.0 44.0)
\clvec(168.33 42.66)(167.33 41.33)(166.0 40.0)
\move(166.0 40.0)
\clvec(164.66 39.0)(163.33 38.33)(162.0 37.66)
\move(162.0 37.66)
\clvec(160.66 37.0)(159.0 36.33)(157.33 36.0)
\move(157.33 36.0)
\clvec(155.66 35.66)(154.0 35.33)(152.33 35.0)
\move(152.33 35.0)
\clvec(150.66 35.0)(149.0 34.66)(147.33 34.66)
\move(147.33 34.66)
\clvec(145.66 34.66)(144.0 35.0)(142.33 34.66)
\move(142.33 34.66)
\clvec(140.66 35.0)(139.0 35.0)(137.33 35.33)
\move(137.33 35.33)
\clvec(135.66 35.66)(134.0 36.0)(132.33 36.33)
\move(132.33 36.33)
\clvec(130.66 36.66)(129.0 37.0)(127.33 37.33)
\move(127.33 37.33)
\clvec(125.66 37.66)(124.0 38.33)(122.33 38.33)
\move(122.33 38.33)
\clvec(121.0 39.0)(119.33 39.33)(118.0 40.0)
\move(118.0 40.0)
\clvec(116.33 40.33)(115.0 41.0)(113.33 41.33)
\move(113.33 41.33)
\clvec(112.0 42.33)(110.33 42.66)(109.0 43.66)
\move(109.0 43.66)
\clvec(107.33 44.0)(106.0 45.0)(104.33 45.66)
\move(104.33 45.66)
\clvec(103.0 46.66)(101.33 47.33)(100.0 48.33)
\move(100.0 48.33)
\clvec(98.66 49.33)(97.33 50.33)(96.0 51.33)
\move(96.0 51.33)
\clvec(94.66 52.33)(93.33 53.33)(92.0 54.66)
\move(92.0 54.66)
\clvec(91.0 56.0)(90.0 57.33)(89.0 58.66)
\move(89.0 58.66)
\clvec(88.33 60.0)(87.66 61.66)(87.33 63.33)
\move(87.33 63.33)
\clvec(87.0 65.0)(87.0 66.66)(87.0 68.33)
\move(87.0 68.33)
\clvec(87.0 70.0)(87.33 71.66)(87.66 73.33)
\move(87.66 73.33)
\clvec(88.0 75.0)(88.66 76.66)(88.66 78.33)
\move(88.66 78.33)
\clvec(89.66 79.66)(89.66 81.33)(90.66 82.66)
\move(90.66 82.66)
\clvec(91.0 84.33)(92.0 85.66)(92.33 87.33)
\move(92.33 87.33)
\clvec(93.33 88.66)(93.66 90.33)(94.66 91.66)
\move(94.66 91.66)
\clvec(95.0 93.33)(96.0 94.66)(96.33 96.33)
\move(96.33 96.33)
\clvec(97.33 97.66)(98.0 99.33)(99.0 100.66)
\move(99.0 100.66)
\clvec(99.33 102.33)(100.33 103.66)(101.0 105.33)
\move(101.0 105.33)
\clvec(102.0 106.66)(102.66 108.0)(103.33 109.66)
\move(103.33 109.66)
\clvec(104.33 111.0)(104.66 112.66)(105.66 114.0)
\move(105.66 114.0)
\clvec(106.33 115.33)(107.0 116.66)(108.0 118.0)
\move(108.0 118.0)
\clvec(108.66 119.66)(109.66 121.0)(110.0 122.66)
\move(110.0 122.66)
\clvec(111.33 124.0)(111.66 125.66)(112.66 127.0)
\move(112.66 127.0)
\clvec(113.33 128.33)(114.33 129.66)(115.0 131.33)
\move(115.0 131.33)
\clvec(116.0 132.66)(116.66 134.33)(117.66 135.66)
\move(117.66 135.66)
\clvec(118.33 137.33)(119.33 138.66)(119.66 140.33)
\move(119.66 140.33)
\clvec(121.0 141.66)(121.33 143.33)(122.66 144.66)
\move(122.66 144.66)
\clvec(123.33 146.33)(124.33 147.66)(125.0 149.33)
\move(125.0 149.33)
\clvec(126.0 150.66)(126.66 152.0)(127.33 153.33)

\move(9 108)
\lvec(181 108)
\move(9 58)
\lvec(181 149)
\move(9 23)
\lvec(150 151)
\move(110 -1)
\lvec(99 151)
\move(155 -1)
\lvec(82 151)
\move(181 27)
\lvec(61 151)

\textref h:C v:C
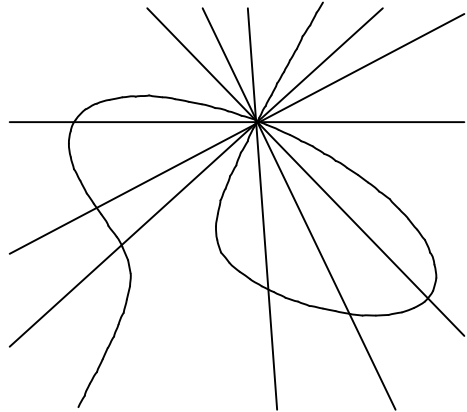
\htext(100 -15){The system $L_{27}(24,6,6,6,6,6,6)$. Each curve is triple.}

}

\end{example}

\begin{example}
We can also use reductions to simplify the computation.
One of the cases treated in
\cite{CCMO}
was the system
$L = \sys_{133}(20^{\times 43})$. In order to compute the dimension
of $L$ using the determinant method, one must consider a $9045 \times 9030$ matrix.
By reduction, it is enough to compute the dimension of the system
$\sys_{D}(20^{\times 5})$, where
$$D = (\overline{36},36,36,36,36,34,33,32,30,27,25,22,18,15,11,7,1).$$
Now the matrix is only $1065 \times 1050$, which allows computation of
the determinant significantly faster.
\end{example}

\section{Special systems with bounded multiplicity}

\begin{theorem}
\label{stillred}
Let $d,M \in \N^{*}$, let $D$ be the diagram obtained by applying a sequence
of $m_{j}$-reductions to the diagram $(\overline{d})$, where $m_{j} \leq M$. If $\# D > 2M(2M - 1)$
then
\begin{enumerate}
\item
$D$ can be written as $(\overline{a},a_{1},\dots,a_{M})$,
$a \geq 2M$,
\item
$D$ is $m$-reducible for any $m \leq M$.
\end{enumerate}
\end{theorem}

\begin{proof}
For the proof (technical but elementary) see
\cite{Dum}.
\end{proof}

\begin{corollary}\label{corstillred}
For a given $m_1,\dots,m_r$ and $d$ large enough, the system $\sys_d(m_1,\dots,m_r)$ is non-special.
\end{corollary}

\begin{proof}Put $M:=\max\{m_1,\dots,m_r\}$, take $d$ such that
$$
\binom{d+2}2>2M(2M-1)+\sum_{j=1}^r\binom{m_j+1}2,
$$
and apply Theorem \ref{stillred}.
\end{proof}

The above theorem allows us to find all special systems with bounded multiplicities
under some conditions.
Fix an $M \in \N^{*}$. The idea is to find a finite set of diagrams $\setofdiag$ with the following
two properties.
\begin{enumerate}
\item
The system $\sys_{d}(m_{1},\dots,m_{r})$, where $m_{i} \leq M$ and
$d$ is large enough, reduces to either system $\sys_{D'}()$ (without
conditions on points) or a system $\sys_{D''}(k_{1},\dots,k_{\ell})$ for
$D'' \in \setofdiag$,
\item
Every system $\sys_{D}(k_{1},\dots,k_{\ell})$, where $D \in \setofdiag$ and
$k_{j} \leq M$, $j=1,\dots,\ell$, is non-special.
\end{enumerate}
The first property can be achieved due to above theorem, which
assures that any diagram large enough can be reduced, and describes
the form of its reduction.
The second property can be checked directly, e.g. by computing
determinants and/or reductions.

Having chosen a suitable set $\setofdiag$, it is enough to look for special
systems in a finite set of systems built on diagrams, which cannot
be reduced to one of the diagrams in $\setofdiag$.

The effectiveness and low time complexity
of computing reductions allowed doing this for multiplicities
bounded by $11$.

The authors wrote and ran suitable programs to find all special systems
up to multiplicity $11$. We present here only the numbers of special
systems (of respective multiplicities). The list of such systems can be found
at
\cite{MYWWW}.
$$
\begin{array}{c|c}
\textrm{ multiplicity } & \textrm{ number of special systems } \\ \hline
\begin{array}{r}
1\\
2\\
3 \\
4 \\
5 \\
6 \\
7 \\
8 \\
9 \\
10 \\
11
\end{array}
&
\begin{array}{r}
0 \\
2 \\
14 \\
91 \\
405 \\
1798 \\
6751 \\
25262 \\
86147 \\
 291868 \\
 929519
\end{array}
\end{array}
$$

As a corollary we have

\begin{corollary}\label{hirsch11}
The Harbourne--Hirschowitz Conjecture holds for systems with multiplicities
bounded by $11$.
\end{corollary}

\section{Bounding the regularity}

\begin{definition}
The \textit{regularity} of a sequence of multiplicities $(m_{1},\dots,m_{r})$
is defined to be
$$\rg(m_{1},\dots,m_{r}) := \min \big\{ d \in \N \ | \ \sys_{d}(m_{1},\dots,m_{r})
\textrm{ is nonempty and non-special} \big\}.$$
\end{definition}

It can be shown (see Corollary \ref{corstillred}) that
regularity is a well-defined natural number. We present here a new algorithm
to find an upper bound on regularity of a system of multiplicities.

\medskip

\begin{algorithm}{\Regularity}

\noindent {\bf Input:} a sequence of multiplicities $(m_{1},\dots,m_{r})$. \\
{\bf Output:} a number $g \in \N$ such that
$\rg(m_{1},\dots,m_{r}) \leq g$. \\

\medskip

\noindent $g \ass -1$ \\
{\bf repeat}

$g \ass g+1$

$d \ass $ \adm $(g, (m_{1},\dots,m_{r}))$\\
{\bf until} $d = \edim \sys_{g}(m_{1},\dots,m_{r}) > -1$ \\
{\bf return} $g$

\end{algorithm}

\begin{theorem}
Algorithm \Regularity stops after a finite number of steps and
returns the upper bound for the regularity of a given system.
\end{theorem}

\begin{proof}
The correctness is obvious. The algorithm must stop in view of
Corollary \ref{corstillred}.
\end{proof}

\begin{definition}
Let $m_{1},\dots,m_{r},k_{1},\dots,k_{r} \in \N$. Define
$$(m_{1}^{\times k_{1}},\dots,m_{r}^{\times k_{r}}) =
(\underbrace{m_{1},\dots,m_{1}}_{k_{1}},\dots,\underbrace{m_{r},\dots,m_{r}}_{k_{r}}).$$
\end{definition}

In the paper \cite{M} a new algorithm for bounding the regularity
is given. We present the comparison of these two algorithms.
We use five sequences of multiplicities (following the author of \cite{M}).
$$\begin{array}{rcl}
L_{1} & = & (500^{\times 2}, 400, 300, 200^{\times 2}, 150^{\times 3},
100^{\times 3}, 80^{\times 3}, 10), \\
L_{2} & = & (350, 300, 250, 200^{\times 2}, 100^{\times 3}, 75^{\times 2},
70, 60, 50^{\times 2}, 30^{\times 2}), \\
L_{3} & = & (386,243,200,170^{\times 2}, 162, 100, 81^{\times 3},
54,40,27^{\times 3}, 25, 10, 9^{\times 3}, \\
& &\phantom{(} 6, 3^{\times 4}, 2, 1^{\times 5}), \\
L_{4} & = & (600, 350, 300^{\times 3}, 180, 150^{\times 4}, 80^{\times 5},
50, 40^{\times 4}, 25, 20^{\times 4}, 16, 10^{\times 5}), \\
L_{5} & = & (500^{\times 3}, 425, 110, 100^{\times 7}, 55, 50^{\times 2},
45, 20^{\times 5}, 10^{\times 4}).
\end{array}$$
For each sequence we give its conjectural regularity
(assuming Hirschowitz--Harbourne Conjecture), and the bounds for
regularity computed by algorithm of Monserrat,
and our reduction method.
Moreover, for each system $L_{i}$, we present also a pair $(d,p)$
(depending on $L_{i}$) such that
$$\dim \sys_{(\overline{d+1},p)}(L_{i}) = \edim \sys_{(\overline{d+1},p)}(L_{i}) > -1$$
(also obtained by reduction method).

$$\begin{array}{|l|c|c|c|c|c|}
\hline
 & L_{1} & L_{2} & L_{3} & L_{4} & L_{5} \\ \hline \hline
\textrm{conjecture} & 999 & 650 & 628 & 959 & 1017 \\
\textrm{\cite{M}} & 999 & 676 & 628 & 985 & 1018 \\
\textrm{reduction} & 999 & 650 & 628 & 959 & 1017 \\
\textrm{reduction} & (998,1) & (649,100) & (627,144) & (958,816) & (1016,833) \\
\hline
\end{array}$$

\section{Reduction and Cremona transformation}

It is worth mentioning that the Cremona transformation can make reducing much easier.

\begin{definition}
Let $L = \sys_{d}(m_{1},\dots,m_{r})$. If $L$ satisfies
\begin{enumerate}
\item $c := m_{1}+m_{2}+m_{3}-d > 0$,
\item $m_{i} \geq c$, for $i=1,2,3$,
\end{enumerate}
then the system $\sys_{d-c}(m_{1}-c,m_{2}-c,m_{3}-c,m_{4},\dots,m_{r})$
is called a \textit{Cremona transformation of $L$}.
\end{definition}

\begin{theorem}
For any system $L=\sys_{d}(m_{1},\dots,m_{r})$, which admits a Cremona
transformation we have:
$L$ is special if and only if it's Cremona transformation is a special system.
\end{theorem}

\begin{proof}
One can show that
if we apply the usual Cremona transformation of $\proj^{2}$
based on points with multiplicities $m_{1}, m_{2}, m_{3}$ then
we obtain a system as above. Of course the dimension cannot change.
\end{proof}

\begin{remark}
A Cremona transformation does not change the expected dimension of
a system (this can be checked by direct computation).
The reductions, however, become smaller and easier to perform.
Many examples have shown that a combination of Cremona transformation
and reduction methods is more effective than using only reductions.
\end{remark}

\begin{example}
Let $L = \sys_{63}(24^{\times 7})$. If we apply the reduction method to
$L$ we obtain that $\dim L \leq 18$. Applying Cremona transformation
to $L$ we get the system from example \ref{ex27}, which gives
$\dim L = 0$.
\end{example}

\end{document}